\documentclass[journal, 10pt]{IEEEtran}

\ifCLASSINFOpdf
\else
\fi

\usepackage{times,amsmath,epsfig}
\usepackage[compress]{cite}
\usepackage{subcaption}
\usepackage{amsfonts}
\usepackage{latexsym}
\usepackage[compress]{cite}
\usepackage{amssymb}
\usepackage{array}
\newcolumntype{C}[1]{>{\centering\arraybackslash}m{#1}}
\usepackage[para,online,flushleft]{threeparttable}
\usepackage{tabulary}
 \usepackage[usenames,dvipsnames]{pstricks}

\begin{document}

\title{A Survey on Energy Trading in Smart Grid}

\author{I.~Safak~Bayram,~\IEEEmembership{Member,~IEEE,} Muhammad Z. Shakir, ~\IEEEmembership{Member,~IEEE,} \newline~Mohamed~Abdallah,~\IEEEmembership{Senior~Member,~IEEE,} and
~Khalid~Qaraqe,~\IEEEmembership{Senior~Member,~IEEE}
\vspace{-25 pt}
\thanks{The authors are with the Department
of Electrical and Computer Engineering, Texas A\&M University at Qatar. Emails:(islam.bayram, muhammad.shakir, mohamed.abdallah, khalid.qaraqe)@qatar.tamu.edu. This publication was made possible by NPRP grant \# 6-149-2-058 from the Qatar National Research Fund (a member of Qatar Foundation). The statements made herein are solely the responsibility of the authors.}}
\maketitle
\begin{abstract}
As the distributed energy generation and storage technologies are becoming economically viable, energy trading is gradually becoming a profit making option for end-users. This trend is further supported by the regulators and the policy makers as it aids the efficiency of power grid operations, reduces power generation cost and the Green House Gas (GHG) emissions. To that end, in this paper we provide an overview of distributed energy trading concepts in smart grid. First, we identify the motivation and the desired outcomes of energy trading framework. Then we present the enabling technologies that are required to generate, store, and communicate with the trading agencies. Finally, we survey on the existing literature and present an array of mathematical frameworks employed.
\end{abstract}


\IEEEpeerreviewmaketitle
\vspace{-10pt}
\section{Introduction}
\subsection{Motivation}
\begin{figure}[t]
\centering
\includegraphics[width=0.75\columnwidth]{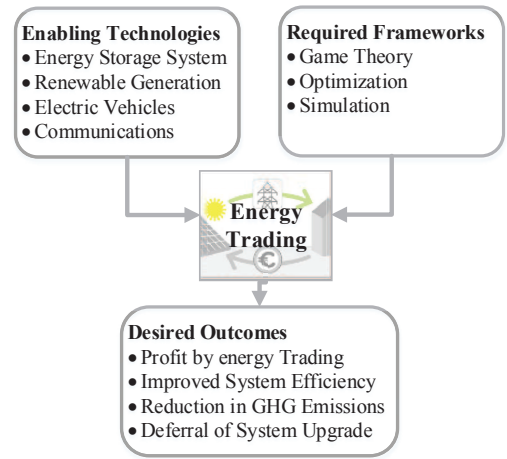}
\caption{An Overview of Energy Trading Components.}\label{chart}
\vspace{-20 pt}
\end{figure}

\vspace{-3pt}
As the world's population is drastically growing, the corresponding additional energy demand (estimated as $35\%$ by year $2040$) should to be continuously supplied in order to sustain the economic development~\cite{exxon}. However, to accommodate the projected demand, the energy efficiency aspect should be carefully considered. The recent advancements in smart grid promise unprecedented improvements in energy efficiency. This can be mainly realized by the deployment of decentralized generation and storage technologies, so that customers can participate in decentralized bilateral energy trading to partially fulfill their needs without further stressing the grid. To that end, in this paper we survey on recent literature on energy exchange and trading in smart grid. An overview of this paper is presented in Figure \ref{chart}.

The success of the aforementioned energy trading framework depends on active participation of end-users and the availability of enabling technologies. The biggest motivation for the users comes in the form of cost savings and profits. From the utility standpoint, the economic benefits are multifaceted \cite{Basu20114348}. The distributed generation, storage, and the exchange of energy will: (1) improve the overall efficiency of the grid operations; (2) minimize system operation cost; and (3) reduce the Green House Gas emissions. Note that the savings in utility operations will further be reflected in electricity tariffs and customers will benefit from lowered tariffs. Next, we explain the aforementioned motivations for energy trading in more details.
\vspace{-5pt}
\subsection{Benefits}

\subsubsection{Improved System Efficiency}

The penetration of distributed generation and storage units and the capability EVs to store considerable amount of energy \cite{speedam}, enable system operators to provide various ancillary services to improve the efficiency of the power grid. For instance, injecting active power from distributed generation options would improve the bus voltages and the power quality, which is mainly realized minimizing voltage sags. Furthermore, with energy trading user demand will be met locally and the use of far-off high capacity generator options will be decreased. This way congestions on transmission lines will be reduced and the corresponding line losses (in the form heat) will be minimized. Corollary, the required system upgrades (due to increasing demand) will be deferred and take place gradually over a wider time horizon. Note that these improvements will enhance the reliability of the equipments and reduce the average customer interruption cost. A detailed survey on improved energy efficiency is presented in \cite{Basu20114348}.
\subsubsection{Reduced System Operation Cost}
\begin{figure}[t]
\centering
\includegraphics[width=0.9\columnwidth]{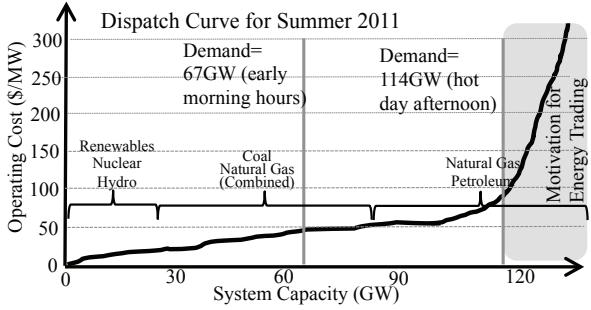}
\caption{Electric Power Generation Cost.}\label{cost}
\vspace{-15 pt}
\end{figure}

In order to meet the stochastic customer demand, utilities dispatch their generation portfolio according their operating cost. Large-scale, low-cost generating units are usually preferred to meet base load demand. On the other hand, as the customer demand increases, system operators dispatch more generators to keep up with the minute-by-minute varying customer demand. During peak hours, that is only around $10\%$ of the day, utilities employ fast-start, high operating cost, and usually gas-power generators to accommodate high electricity demand. One of the main motivations of the energy trading is to reduce to peak-to-average demand ratio by locally trading energy during peak hours. In Figure \ref{cost}, we present a typical cost of dispatch curve for summer 2011 in the US~\cite{eiaCurve}. During the peak hours, the system operation cost increases exponentially. If the energy trading takes places in the shaded region, both the utility and the users will benefit in the form of cost savings.
\subsubsection{Reduction in GHG Emissions}

Countries around the globe have set targets for reducing the GHG emissions in the next decade. For instance, European Union countries are targeting to cut the GHG emission by $20\%$ by year $2020$. The full exploitation of renewable and distributed energy generation is a critical element in reducing the GHG emissions since the power grid operations are responsible for the one fourth of the global GHG emissions. This can be mainly realized by the selling/buying excess generation at micro-grid level during the peak hours and improving the overall system efficiency. On-site power generation in the meantime often releases significant amounts of waste heat, which can be recycled for heating and cooling buildings, refrigeration through absorption chilling, and heating water. Such utilization can improve overall "energy efficiencies" of consumer from around 35\% to as high as 85\%, with additional reductions in per capita CO2 emissions \cite{shakir1}. It is shown in \cite{pratt2010smart} that the aforementioned smart grid features can reduce the carbon emission by $12\%$ by year $2030$. A detailed analysis on the GHG savings with the smart grid features is presented in \cite{grid2008green}.
\vspace{-10pt}
\subsubsection{Energy Profiling}
\vspace{-3pt}
The efficiency of micro-grid generators is greatest when they are networked together with smart virtual micro-grids to profile the user demand and subsequently control the flow and generation of the energy.  Like the bulk power grid, smart micro-grids generate, distribute, and regulate the flow of electricity energy to consumers, but do so locally. Being local they offer much lower line losses in comparison with the higher line losses associated with the energy transmission over long-distance.  Being networked is an ideal way to integrate renewable resources on the community level and allow for customer participation in the electricity enterprise, balancing supply and demand near point of use. The emerging energy trading processes form the building blocks of the perfect power system of future. Perhaps one of the greatest advantages of networked micro-generators is that they can be better adapted to meet the needs of the future energy demand and CO2 reduction targets, and provide rapid response to balancing between demand and supply at smaller granularity than the central grid can offer. Rather than relying on utility companies to build capital-intensive, full-scale power plants, networked micro-generators can enlarge the overall electricity supply rapidly and cost-effectively using relatively small local generators.
\section{Enabling Technologies \& Models}
\vspace{-5pt}
\subsection{Distributed Renewable Generation}
Micro power stations (micro-grid generation) are modern, small, on-site distributed energy generators that can operate grid-connected or be isolated from it. They generally have capacities under 10 megawatts (MW) using renewable energy sources, such as solar panels and wind turbines, or high efficiency conversion of bio-energy or fossil fuels. Over the past decade there has been a strong push to accelerate the integration of distributed renewable energy sources into the existing generation portfolio. With suitable control they are able to reduce peak loads and can provide reliable power for commercial buildings, industrial facilities, residential neighborhoods, college/university campuses, and military bases. What is more, they don't require connection to a national grid. Hence, they are attractive technologies in achieving specific local operational objectives, such as reliability, carbon emission reduction, diversification of energy sources, and cost reduction, established by the community being served.

\begin{figure}[t]
\centering
\includegraphics[width=0.9\columnwidth]{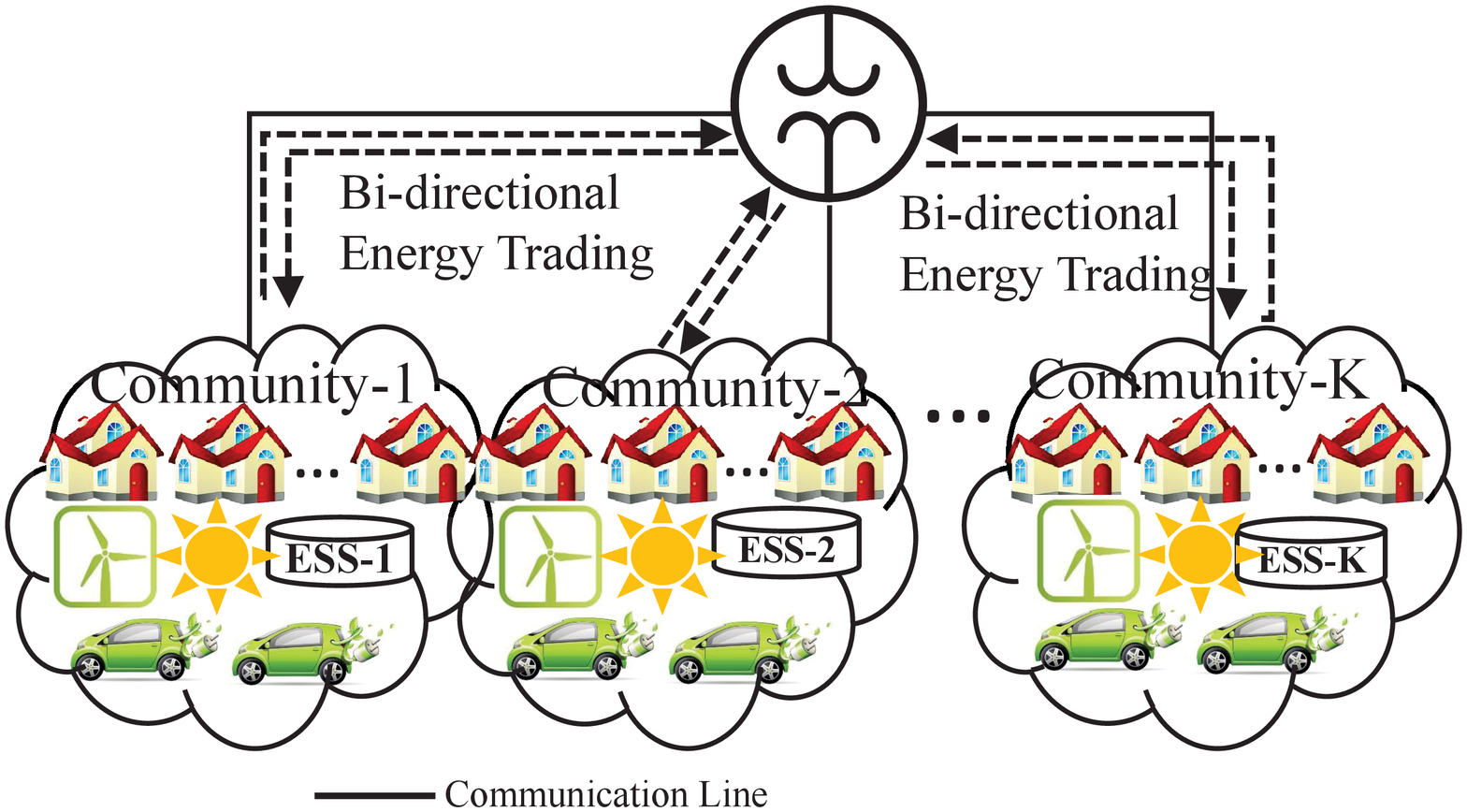}
\caption{Energy Trading Among Microgrids}\label{systemModel}
\vspace{-20 pt}
\end{figure}

To this end, microgrid architectures are key enablers in reaching these goals as they employ renewable generation and offer flexible energy management solutions. Microgrid users, either individuals or groups, act like ``prosumers'' who can produce and consume energy at the same time. Moreover during the periods of supply-demand mismatch, ``prosumers'' can interact among each other and trade electric generation over a marketplace (an overview is illustrated in Figure \ref{systemModel}). On one hand, users make extra profit by selling their excess power or buying cheap electricity from their neighbors. On the other hand, the elimination of starting fast ramp generators reduces the operation cost of the power system. Hence,  energy trading creates a win-win situation for both parties. In literature, in order to quantify the additional savings, second order functions are predominantly used to model the power generation cost of such generators, which is given by:
\begin{equation}\label{genCost}
C(x) = {a_1}{x^2} + {a_2}x + {a_3},
\end{equation}
where $x$ is the active power output and $a_{1}$, $a_{2}$, and $a_{3}$ are the cost coefficients of the generator. \eqref{genCost} usually serve as a part of the profit function in energy trading literature \cite{TusharTSG,Tushar}.
Another important aspect of the research efforts in renewable generation is to choose appropriate stochastic models to capture the intermittency in energy trading. The power generation in wind turbines are correlated with the wind speed which is modeled by a Weibull distribution, and the scale parameters are associated with strength of the wind and the peak of the wind distribution. A similar mathematical model to predict the output of solar power generation is presented in \cite{liangstochastic}, where the output power is computed by the following three parameters; efficiency of the photovoltaic (PV) system, solar cell temperature, and the intensity of the solar radiation.

On the other hand, from system design point of view, models in stochastic network calculus are successfully integrated in smart grid application for the renewable energy generation\cite{sLow}. In this frameworks, renewable generation is modeled as a general stochastic arrival process and it is stored in an energy storage. Then the goal is provide a probabilistic bound on the percentage of time the supply will be short of demand. Corollary, there has been a growing body of literature on the use of renewable generation in energy trading frameworks \cite{Vytelingum,ilic,Shroff,6688037,6798766,TusharTSG,6477198,6485958}.
\vspace{-5pt}
\subsection{Energy Storage Systems (ESS)}
As the ESS technology is becoming more economically viable, the role of ESS in energy trading will be more prominent. For large-scale renewable generation (e.g., solar arrays, wind farms), the ESS will be used to smooth out the output of the system \cite{teleke1}. On the other hand, for end-user applications (distributed) community-based energy storage systems have already gained popularity~\cite{roleOfESS}. In this case, the goal is to deploy small size storage units in the residential feeders to accommodate the demand of several houses during peak demand. Similarly for energy trading applications, the primary role of the energy storage system will be to store off-peak hour renewable energy, so that users can use and exchange it during the periods of peak demand. The bibliography on energy trading applications with ESS include \cite{Vytelingum, ilic,Shroff,6688037,6798766,6060940,Tushar,6477198}. Overall the goals of the ESS technologies are: (1) improve power grid optimization for bulk power production; (2) balance the power system operations with intermittent and/or diurnal renewable generation options; (3) help to defer capital-intensive upgrades in the transmission and distribution grids; (4) provide ancillary services to grid operations.

From the modeling perspective researchers threat storage units as linear entities and use resource provisioning literature from communications to solve the sizing problems~\cite{sgc11,sgc12,jsac,sgc13,springer}. To that end, the role of ESS modeling in optimization problems depends on the underlying assumption. If it is assumed that the ESS has already been acquired and operated by utility, then the energy trading entities use the ESS size as a constraint in the optimization problem. However, if the owner of the ESS is also the energy trading entity (micro grid or individual end-users), then the size of the ESS becomes an additional cost term in the objective function.
\vspace{-5pt}
\subsection{Electric Vehicles}
Even though the primary goal of Electric Vehicles (both pure electric and plug-in hybrid) is to offer environmentally friendly and cost-effective transportation options, the capability of EVs to store huge amount of electric power (e.g., the size of Nissan Leaf battery can store up to two household demands in the U.S.) makes them a natural player in energy exchange mechanism. With the use of bidirectional chargers, EVs can exchange electric power with the power grid or other EVs\cite{V2VProc}. From energy trading standpoint, there are three emerging concepts on the use of EVs. The first one is Vehicle-to-Home (V2H), in which the vehicle battery pack acts exactly the same as the stand-alone ESS given in the previous section. As a second scenario, EVs can exchange energy among each other (Vehicle-to-Vehicle, V2V) to fulfill their requirements in a cost effective way. The most popular EV application concept is the Vehicle-to-Grid (V2G) where the stored energy is exchange with the grid. The predominant use of EVs in the literature is for energy trading to make extra profit by selling excess power and the related literature includes \cite{Schaar,6102331,6798766,TusharTSG,6193525,bayramproviding}, however the stored energy of EVs can also be used in ancillary services \cite{6060940,6084772,DBLP:journals/tsg/SortommeE12a}.

In energy trading applications, the optimization problems may include different objectives. The most profound ones are: (1) cost minimization; (2) cost-emission minimization; (3) power-loss minimization; and (4) peak-to-average load minimization. Depending on the objective, different constraints may be required for the energy trading schemes. The most important constraint is of course the electricity prices which comes from the utility company. In order to minimize the power losses, the physical distance between the trading entities needs to be considered. Also technological constraints like State-of-Charge limitations, battery types, battery ratings and efficiency will play a critical role.
\vspace{-5pt}
\subsection{Two-way Communications}
The success of the energy trading mechanism heavily depends on the availability of the necessary communication infrastructures to ensure reliable information dissemination. In energy trading, participants need to update their demand or the amount of available energy to sell with the market place via two-way communication technologies. Also communication networks will enable trading entities to monitor their renewable generation and the state of charge at the storage unit. The literature in its current state assumes that there is perfect communication between all players. However, it is also important to quantify the impacts of communication system performance on the operations of the energy trading mechanisms as it will create another level of uncertainty. Trading entities can make sub-optimal choices due to loss of communications.	 For instance, the work presented in \cite{6231158} quantifies the effects of communication system reliability in power consumption of smart grid users. Similar concepts can be used in energy trading programs.
\begin{figure}[t]
 \begin{threeparttable}
\centering
 \caption{Literature Review}\label{table1}
\begin{tabular}{| C{1.55cm}| C{0.8cm}| C{0.45cm} | C{0.45cm} | C{3.45cm}   | }
\hline
	Ref. &Rnwbl. &ESS&EVs &Model  \\ \hline
	\cite{Vytelingum,ilic}&Y & N& N & Double Auction  \\ \hline
\cite{Shroff,6688037}&	Y & Y& N & Stochastic Optimization  \\ \hline
\cite{Schaar,6102331}&	N& N & Y& Noncooperative Game  \\ \hline
\cite{cui2014electricity}&	Y & N &N & Social Welfare Max. \\ \hline
\cite{6798766}&	N &Y & Y & Double Auction \\ \hline
\cite{6060940}&	N & N & Y &  Noncooperative Game \\ \hline
\cite{Tushar}&	Y & Y& Y & Stackelberg Game  \\ \hline
\cite{TusharTSG}&	N & N & Y & Stackelberg Game \\ \hline
\cite{6193525}&	N & N & Y & Double Auction  \\ \hline
\cite{6477198} &	Y &Y & Y & Bidding  \\ \hline
\cite{6485958} &	Y &N &N & Convex Optimization  \\ \hline
\cite{immune} &	Y &Y & N & Particle Swarm Optim.  \\ \hline
\end{tabular}
\begin{tablenotes}
\scriptsize
\item Abbreviations: References (Refs.), Renewable Generation (Rnwbl.), Energy Storage Systems (ESS), Electric Vehicles (EVs)
\end{tablenotes}
\end{threeparttable}
\vspace{-20pt}
\end{figure}
\section{Required Frameworks}
The literature on energy trading can be classified into several subcategories by considering the different combinations of employed enabling technologies that are presented in the previous section. Another important aspect in categorizing the literature is the employed mathematical framework. In general, such frameworks can be classified into three categories. If the energy trading scenario is set to have only one user or a central controller who can dictate his decisions to a group of users, then the appropriate framework would be to use single objective maximization tools such as convex, stochastic, or swarm optimization, or social welfare maximization. However, in majority of the cases there are multiple interacting users who try to optimize their own utilities without considering the rest of the user and the grid conditions. In this case game theoretic approaches are adopted to find the optimal solutions in a decentralized way. As a final approach, we present the literature on simulation-based solutions. To that end, we extent our categorization by adding another dimension present the most popular approaches used in the existing literature and the summary of the literature is presented in Table~\ref{table1}.
\vspace{-5pt}
\subsection{Decentralized Solution: Game Theoretic Approach}
\subsubsection{Auction Mechanism}
Auction mechanisms have been the cornerstone of many applications in wholesale and retail electric power markets. Similar to traditional auction mechanism, the primary goal of the distributed energy trading is to find the lowest-cost matching between the supply and demand to maximize the economic efficiency \cite{morey2001power}. In \cite{6084772}, authors present a generic auction mechanism for energy trading in local markets. Moreover, the works presented in \cite{Vytelingum,ilic,6798766,6193525} employ double auction mechanism in their model. More specifically, the work presented in \cite{Vytelingum} proposes a market mechanism using continuous double auction aiming to improve the market efficiency. In this framework buyers (bid) and sellers (ask) continue trading period until the market is cleared that is when a bid exceeds an ask. Similarly, authors of \cite{6798766} model the energy trading among distributed energy storage units with game theory using double auction mechanism. Moreover, \cite{6193525} formulates a game for energy exchange among electric vehicles and the power grid with double auction mechanism.
\subsubsection{Stackelberg Game}
In economics Stackelberg game models the behaviour of two agents, one of them being the \emph{leader} with the first move advantage and the other one being the \emph{follower} who plays a best response strategy to maximize his own utility \cite{globecom}. In energy trading applications, the aggregator often becomes the \emph{leader} and sets the prices according to the needs of the market and aims to motivate users for participation. For instance, the work presented in \cite{Tushar} incentives customers to sell their surplus energy during the peak hours. Moreover, authors of \cite{TusharTSG} uses \emph{leader-follower} game to model the energy exchange in vehicle to grid application and the solution of the game is proven to be the socially optimal point. In both papers, utility functions are used to capture the behaviour of players, and the strategies are developed to maximize this functions which are concave.

\subsubsection{Non-cooperative Game}
Non-cooperative game theory has been used extensively in many applications in economics and engineering to study the interaction among independent and self-interested agents. In energy trading applications the non-cooperative games are employed to calculate the amount of energy to be sold in the market and the optimal solution (if exists) is the Nash equilibrium, where no player has incentive to deviate from his strategy. The work presented in \cite{Schaar} uses non-cooperative game to solve the optimal amount of energy exchange among a group of plug-in hybrid vehicles. Similar approach has been adopted in \cite{6102331} and \cite{6060940}.
\subsection{Centralized Solution: Single Objective Maximization}
In the case of centralized approaches, the trading agencies act as one entity or follow the orders of a central controller who is assumed to know all the information about buyers and sellers. In this case, single objective maximization techniques are adequate to compute the optimal amount of energy to be traded. In \cite{6485958} authors model the peer-to-peer energy trading in a microgrid environment and use convex optimization tools to minimize the total energy generation and transportation cost. Similarly the authors in \cite{immune} employ particle swarm optimization schemes to minimize the fuel cost generation cost. On the other hand, in the case of renewable generation stochastic optimization techniques have been used to address the uncertainty in generation. The work presented in \cite{Shroff} proposes a profit maximization problem from end-users standpoint using stochastic programming.
\subsection{Simulation-based Solutions}
The third group of approaches use simulation based studies to model the behavior of multi scale decision making agents. The main part of such approaches is the use of statistical learning algorithms (e.g., reinforcement learning, Q-learning etc.), so that trading agents can derive long-term profit making policies in an autonomous way \cite{Reddy:2011:SLA:2283516.2283637,reddy2011learned,peters2012autonomous,vytelingum2010agent}. In \cite{peters2012autonomous} authors use electricity brokers as aggregators to manage the balance between buyers and sellers. A similar approach is used in \cite{reddy2011learned} and the broker agents behavior is modeled with Markov Decision Processes and Q-learning techniques. Furthermore, the work presented in \cite{japan} proposes a simulation based modeling for local energy trading.
\vspace{-5pt}
\section{Conclusions}
Over the past few years, there has been a growing interest in energy trading applications in smart grid. In this paper, we surveyed on the distributed energy trading and exchange in smart grid. We characterized the enabling technologies as the renewable generation, energy storage, electric vehicles, and communications systems. Furthermore, we divide the mathematical frameworks into three. The first group included game theoretic models that are used for multi-agent decision making. On the other hand, second type of approaches uses single objective maximization. Finally, we showed the use of the simulation-based studies for energy trading mechanism.
\ifCLASSOPTIONcaptionsoff
  \newpage
\fi

\bibliographystyle{IEEEtran}
\bibliography{globalSIP}

\end{document}